\documentclass{article}
\usepackage{amsfonts,amssymb,latexsym}
\newcounter{num}[section]
\newcommand{\Num}{\refstepcounter{num}%
\textbf{\arabic{section}.\arabic{num}}}

\newcommand{\Theorem}{\textbf{Theorem~}}
\newcommand{\Proof}{{\emph{Proof}}}
\newcommand{\Def}{\textbf{Definition~}}
\newcommand{\Conj}{\textbf{Conjecture~}}
\newcommand{\Lemma}{ \textbf{Lemma~}}

\newcommand{\Prop}{\textbf{Proposition~}}

\newcommand{\Notation}{\textbf{Notations~}}

\newcommand{\Kc}{{\cal K}}
\newcommand{\Xc}{{\cal X}}
\newcommand{\Ac}{{\cal A}}
\newcommand{\Bc}{{\cal B}}

\newcommand{\Ch}{{{\mathfrak C}{\mathfrak h}}}

\newcommand{\al}{{\alpha}}
\newcommand{\la}{{\lambda}}
\newcommand{\De}{{\Delta}}

\newcommand{\eps}{{\varepsilon}}

\newcommand{\Tn}{{\mathrm{T}}(n,\Fq)}

\newcommand{\GL}{{\mathrm{GL}}}
\newcommand{\Irr}{{\mathrm{Irr}}}
\newcommand{\Ind}{{\mathrm{Ind}}}
\newcommand{\Ad}{{\mathrm{Ad}}}
\newcommand{\row}{{\mathrm{row}}}
\newcommand{\col}{{\mathrm{col}}}

\newcommand{\rt}{{\mathrm{right}}}
\newcommand{\lt}{{\mathrm{left}}}

\newcommand{\JDDr}{{J_{\Db,\rt}}}
\newcommand{\NDDr}{{N_{\Db,\rt}}}
\newcommand{\RDDo}{{R_{\Db}^\circ}}
\newcommand{\RDDr}{{ R_{\Db,\rt}}}
\newcommand{\RDDl}{{ R_{\Db,\lt}}}
\newcommand{\laDD}{{ \la_{\Db}}}
\newcommand{\JDDrkl}{{J_{\Db,\rt}^{k,\ell}}}

\newcommand{\tP}{\tilde{P}}

\newcommand{\Fq}{{\Bbb F}_q}
\newcommand{\Cb}{{\Bbb C}}

\newcommand{\Db}{{\Bbb D}}

\renewcommand{\leq}{\leqslant}
\renewcommand{\geq}{\geqslant}

\begin{document}
\Large

\title{Towards a supercharacter theory of the parabolic subgroups }
\author{A.N.Panov
\thanks{The research is supported by the grant RSF-DFG 16-41-1013}\\
\\
 Department of Mathematics, 
Samara  University,\\
ul.~Akademika Pavlova 1, Samara, 443011, Russia\\
email: apanov@list.ru}
\date{}

 \maketitle

\begin{abstract}
	
A supercharacter theory is constructed for the parabolic subgroups of the group  $\mathrm{GL}(n,\Fq)$ with blocks of orders less or equal to two. The author formulated  the hypotheses on construction of  a supercharacter theory for an arbitrary parabolic subgroup in  $\mathrm{GL}(n,\Fq)$.

\end{abstract}

\section{Introduction}

The problem of classification of irreducible characters is an extremely difficult, "wild"\'~ problem for some groups like  the unitriangular $\mathrm{UT}(n,\Fq)$ and the parabolic subgroups in  Chevalley groups. For these groups, it seems logical to construct a first approximation of the theory of irreducible characters, which is called a supercharacter theory.

The notion of a supercharacter theory was suggested by P.Diaconis and  I.M.Isaacs in 2008 in the paper   \cite{DI}. A priori a group admits several supercharacter theories, one of the ones is the theory of irreducible representations. If it is impossible to   classify the irreducible charcaters for a given group, then  the goal is to construct a supercharacter theory that provides the best approximation of the theory of irreducible characters. One of  examples of such a theory  is the theory of basic characters of
   C.Andr\'{e} (see \cite{A1,A2,A3,A4,Yan}) for the unitriangular group and the supercharacter theory for the algebra groups   \cite{DI} (by definition, an algebra group is a group of the form  $1+J$, where $J$ is a finite dimensional associative nilpotent algebra).

    Many papers were devoted to construction of supercharacter theories for different groups and to development of  the general theory. Observe a few of them:
  supercharacters for abelian groups and their application in the number theory  \cite{Number-1, Number-2},
the superinduction for the algebra groups \cite{ThiemB, ThiemRest, MT}, the super\-cha\-rac\-ter theory for Sylov subgroups  of the orthogonal and symplectic groups over a finite field \cite{AN}, application in the random walk problem on groups \cite{Walk}, the supercharacter theory  for semidirect products \cite{H,SA}, characterization of the Hopf algebra of supercharacters for the unitriangular group  \cite{VERY}. One can found the bibliography in the paper  \cite{VERY}.

In the most of these papers the supecharacter theories were constructed for  unipotent group  closely related to the unitriangular group. One of the main  goals is to enlarge the list groups that admit
a good supercharacter theory.
From the unitriangular group  it is natural to move to the construction  of supercharacter theories for the Borel and the parabolic subgroups  in  $\GL(n,\Fq)$.
The main problem is to construct a supercharacter theory, which would be a natural continuation of the supercharacter theory for algebra groups by P.Diaconis and  I.M.Isaacs.

In the series of papers  \cite{P1,P2,P3,P4},  the author constructed the supercha\-rac\-ter theory for the finite groups of triangular type; a  special case of these groups is  the triangular group  $\mathrm{T}(n,\Fq)$. This theory provides the better approximation of the theory of irreducible representations than the general construction of supercharacters for semidirect products from the paper  \cite{H}.
In the above papers, we calculate supercharacter values on superclasses for the triangular group \cite{P1}, research the restriction and induction in the framework of  constructed supercharacter theory, prove the  Frobenius reciprocity theorem for supercharacters \cite{P3}, obtain an analog of A.A.Kirillov formula for supercharacters \cite{P2}, characterize  the  Hopf algebra of supercharacters of the triangular group as partially symmetric functions in noncommuting variables \cite{P4}.

In the present paper, we formulate  conjectures on  a supercharacter theory for the parabolic subgroups of the group  $\GL(n,\Fq)$ (see Conjectures  \ref{conjone}, \ref{conjtwo}, \ref{conjthree}). We prove these conjectures  for the parabolic subgroups with blocks of orders  $\leq 2$ (see Theorems  \ref{classlow}, \ref{classlowstar},  \ref{theosupercl}, \ref{maintheo}).
Some statements are proved for an  arbitrary parabolic subgroup
(see Propositions  \ref{blockJ} and \ref{sss}).
 In the constructed supercharacter theory, we combine the theory of irreducible characters  of the group  $\GL(2,\Fq)$ and the theory of basic characters of  C.Andr\'{e}.

\section{Notion of a supercharacter theory}

Recall the notion of a supercharacter theory from the paper   \cite{DI}. Let  $G$ be an arbitrary finite group, ~ $1\in G$ be the unit element. Let we have a pair $(\Ch, \Kc)$, where  $\Ch = \{\chi_1, \ldots, \chi_m\}$  is a system of pairwise disjoint characters of the group $G$, and  $\Kc = \{K_1,\ldots, K_m\}$ is a partition of the group   $G$.   \\
\\
 \Def\Num\label{defsuptheory} The pair $(\Ch, \Kc)$ is said to determine  a supercharacter theory if it satisfies the following conditions:\\
1) ~ each character  $\chi_i$ is constant on each subset  $K_j$;\\
2) ~ $\{1\} \in \Kc$.
\\

The characters of    $\Ch$ are called    {\it supercharacters}, and the subsets from     $\Kc$ are    {\it superclasses}.

Observe that the number of superclasses equals to the number of supercharacters. The square table of supercharacter values on superclasses is called a \emph{supercharacter table}.

Let  $X_i$ denote the \emph{support of character} $\chi_i$ (i.e. the set of all irreducible components of  $\chi_i$).   It can be seen the  importance of condition  2)  from the following lemma.\\
\\
\Lemma\Num~ \cite[Lemma 2.1]{DI} \textit{Let the system of disjoint characters   $\Ch$ and the partition  $\Kc$ obey the condition  1).
Then the condition   2) is equivalent to the following condition\\
2')~ the system of supports   $\Xc=\{X_1,\ldots, X_m\}$ is a  partition of the set of irreducible characters $\Irr(G)$; moreover, each character $\chi_i$ equals to  the character  $$\sigma_i = \sum_{\psi\in X_i} \psi(1)\psi $$
up to a constant factor.}\\

So, the supercharacter character  defines a pair of partitions  $(\Xc,\Kc)$ (where $\Xc$ is a partition of    $\Irr(G)$, and  $\Kc$ is the one of the group  $G$)
 of equal number of components such that  each character  $\sigma_i$ is constant on each $K_j$.
One can take this property as a definition of a supercharacter theory.

Observe some other properties of the systems    $\Ch$ and  $\Kc$.\\
\\
 \Prop\Num ~\cite[Theorem  2.2.]{DI} \\
\textit{1) Each superclass is a union of conjugacy  classes.\\
2) The partition  $\Kc$ is uniquely defined by the partition $\Xc$ and vice versa.\\
3) The principal character is a supercharacter (up to a constant factor).}

\section{Conjectures for the parabolic subgroups of  $\GL(n,\Fq)$}\label{hyp}

 Let $P$ be a parabolic subgroup of  $\GL(n,\Fq)$ with blocks of orders $n_1,\ldots, n_s$ that contains the group of upper triangular matrices $\Tn$.
The subgroup  $P$ is a semidirect product   $P=RN$, where  $R$ is a direct product of the subgroups  $R_i=\GL(n_i,\Fq)$, ~ $i=1,2,\ldots,s$, and  $N=1+J$ is an unipotent subgroup. Observe that $J$ is an associative subalgebra invariant with respect to  the left and right multiplication of  $P$. The subgroup   $N$ is an algebra group;
 it admits  the supercharacter theory by  P.Diaconis and  I.M.Isaacs.

We define a \emph{root} as an arbitrary pair of positive integers  $(i,j)$, where $1\leq i,j\leq n$, ~ $i\ne j$. For the root $\al=(i,j)$, the number  $i$ will is called a \emph{number of row of the root} $\al$ (notation $i=\row(\al))$. Respectively,  $j=\col(\al)$ is the \emph{number of column of the root}  $\al$.

 The Weyl group $W=S_n$ naturally acts on the set of roots  by   $w(i,j)=(w(i),w(j))$.
 The root $(i,j)$ is  \emph{positive} if  $i<j$. For each root $\al=(i,j) $, we denote by  $E_\al$  the corresponding matrix unit.
Consider the subset  $\De_J$ of positive roots such that  $\{E_\al:~ \al\in \De_J\}$  is a basis of  $J$.
The action of    $W$ on roots induces the action of the Weyl group  $W_R$ of the subgroup  $R$ on $\De_J$.

 Introduce the following equivalence relation on  $J$.\\
 \\
\Def\Num~  The elements  $x, x'\in J$  are \emph{equivalent} if there exist the elements  $r\in R$ and $a,b\in N$ such that
$x'=raxbr^{-1}$.

Observe that the equivalence classes coincide with the orbits of the group  $\tP=R\ltimes(N\times N)$ on $J$.\\
\\
\Def\Num\label{basic}~
A subset  $D\in \De_J $ is a  \emph{ rook placement} if  there is  no more then one root of $D$
in each row and column ( C.Andr\'{e} used the  term  a "basic subset").

For the rook placement    $D$,  we construct the element
\begin{equation}\label{xD}
x_D=\sum_{\al\in D} E_\al.
\end{equation}
\Conj\Num\label{conjone} \\
1) For any  $x\in J$ there exists a  rook  placement  $D$ in $\De_J$ such that   $x$ is equivalent to   $x_D$.\\
2) Two elements  $x_D$ and  $x_{D'}$ are equivalent  if and only  if  $D$ and  $D'$ are conjugated  with respect to the action of  $W_R$ on $\De_J$.\\

Define the left and right actions of the group  $P$ on the dual space  $J^*$ by the formulas
$g\la(x)=\la(xg)$ and $\la g(x)=\la(gx)$. This actions enable us to define an equivalence relation on  $J^*$.
\\
\\
\Def\Num~  We say that the elements  $\la, \la'\in J^*$  are \emph{equivalent } if there  exist the elements  $r\in R$ and $a,b\in N$ such that
$\la'=ra\la br^{-1}$.\\

As above the equivalence classes coincide with the orbits of action of the group  $\tP$ on $J^*$.
 Let $\{E^*_\al:~ \al\in\De_J\}$ be a dual basis with respect to the basis  $\{E_\al\}$ in $J$.
 For the rook placement   $D$, we define the element
\begin{equation}\label{laD}
\la_D=\sum_{\al\in D} E^*_\al.\end{equation}
Let us formulate the conjecture dual to Conjecture  \ref{conjone}.\\
\\
\Conj\Num\label{conjtwo} \\
1) For any  $\la\in J^*$ there exists a  rook placement  $D$ in $\De_J$ such that   $\la$ is equivalent to   $\la_D$.\\
2) Two elements  $\la_D$ and  $\la_{D'}$ are equivalent  if and only if  $D$ and  $D'$ are conjugated  with respect to the action of  $W_R$ on $\De_J$.\\

Let us define an equivalence relation on the group  $P$.\\
\\
\Def\Num\label{superclass} We say that the elements  $g, g'\in P$ are \emph{equivalent} if there exist  the elements  $r\in R$ and $a,b\in N$ such that
$$g' = 1+ra(g-1)br^{-1}.$$
Each equivalence class  $K(g)$ is an orbit of the group  $\tilde{P}$ on $P$.\\

Let us construct  the system of characters. Let  $D$ be a rook placement in $\De_J$, and  $\la_D$ be the corresponding linear form.
Consider the stabilizer  $N_{D,\rt}$ of  $\la_D$ with respect to the right action on the group  $N$ on $J^*$.
 The subgroup $N_{D,\rt}$ is an algebra group  $N_{D,\rt} = 1+ J_{D,\rt}$, where $J_{D,\rt}$ consists of all  $x\in J$ such that  $\la_D(xy)=0$ for any  $y\in J$.

  Let $R_{D,\rt}$ (respectively,  $R_{D,\lt}$) be a stabilizer of  $\la_D$ with respect to the right (respectively, left) action of the group  $R$ on $J^*$. Denote by  $R_D^\circ$ the subgroup
$$ R_D^\circ = R_{D,\rt} \cap R_{D,\lt}.$$
Form the subgroup  $P_D=R_D^\circ N_{D,\rt}$. Each element of this subgroup can be uniquely  presented in the form $g=r+x$, where  $r\in R_D^\circ$ and $x\in J_{D,\rt}$.

Fix a nontrivial character   $t\to \eps^t$ of the additive group of the field  $\Fq$ with values in the multiplicative  group  $\Cb^*$.
Let  $T$ be a representation of the subgroup  $R_D^\circ$ with the character  $\theta$. Denote by  $\Xi$
 the representation of the group  $P_D$ defined as follows
 $$\Xi(g)=T(r)\eps^{\la_D(x)},$$
 where $g=r+x$,~ $r\in R_D^\circ$,~ $x\in J_{D,\rt}$.
 Let us show that $\Xi$ is really a representation:
$$\Xi(gg') = \Xi((r+x)(r'+x')) = \Xi(rr' + r'x + x'r  + xx') =$$
$$T(rr')\varepsilon^{\la_D(r'x)} \varepsilon^{\la_D(x'r)}\varepsilon^{\la_D(xx')} = T(r)T(r')\varepsilon^{\la_D(x)} \varepsilon^{\la_D(x')} = \Xi(g)\Xi(g').$$
 The representation $\Xi$ has  the character
$$
    \xi_{D,\theta}(g) = \theta(r)\varepsilon^{\la_D(x)}.
$$
Consider the induced character
\begin{equation}\label{superch}
\chi_{D,\theta} = \Ind(\xi_{D,\theta}, P_D, P).
\end{equation}
The subgroup  $R_D^\circ$ is a normal subgroup in  $$ R_D=\{r\in R:~ r\la_Dr^{-1}=\la_D\}.$$
A character  $\theta$ of the subgroup  $R_D^\circ$ is called  $R_D$\emph{-irreducible} if  $$\theta=\sum_{r\in R_D/R_D^\circ} r\psi r^{-1}$$
for some irreducible character $\psi$ of the subgroup   $R_D^\circ$.
\\
\\
\Conj\Num\label{conjthree}~ The system of characters  $\{\chi_{{D,\theta}}\}$ (here  $D$ runs through the set of representatives of $W_R$-orbits on rook placements in $\De_J$, and  $\theta$ runs through the set of $R_D$-irreducible characters of the subgroup  $R_D^\circ$) and the partition of  $P$ into
the equivalence classes  $\{K(g)\}$ (see definition  \ref{superclass}) give rise to a supercharacter theory of the parabolic subgroup  $P$.

\section{The equivalence classes in  parabolic subgroups}

Let  $P$  be a parabolic subgroup in $\GL(n,\Fq)$ of type  $(n_1,\ldots, n_s)$.
Divide the segment  $I=[1,n]$ into the consecutive segments $I=I_1\cup\ldots \cup I_s$ of sizes $n_1,\ldots, n_s$. The set of roots  $\De_J$ is divided into the blocks  $$\De_J=\bigcup_{1\leq k<\ell\leq s} I_k\times I_\ell.$$
For each  $1\leq k\leq s$, we define a  $k$th \emph{block-row } in $\De_J$ as  the union of blocks $I_k\times I_j$ over all $k<j\leq s$. Respectively, a $\ell$th \emph{block-column} is the union  $I_i\times I_\ell$ over all  $1\leq i <\ell$.
\\
\\
\Def\Num\label{assblock} \\
1) Let  $A_k\subseteq I_k$ and $A_\ell \subseteq I_\ell$  for $k<\ell$.
If  $|A_k|=|A_\ell|$, then we refer to the subset  $A_k\times A_\ell$  in $I_k\times I_\ell$ as a  \emph{block-rook}.\\
2) A subset $\Db$ in  $\De_J$ is a  \emph{ block-rook placement} if
$\Db$ is a union of block-rooks, which do not attack each other.
The last means that each row and each column intersects with no more than one block-rook. \\
3) We say that the  element  $x_\Db$ in $J$ is  \emph{associated with a block-rook placement} $\Db$ if
$$ x_\Db = \sum_{(i,j)\in\Db} x_{ij}E_{ij},$$
 moreover, for each block-rook  $A_k\times A_\ell\subset \Db$ the submatrix  $$\{x_{ij}:~ i\in A_k,~ j\in A_\ell\}$$ is non-degenerate.\\
4) We say that the  element $\la_\Db$ in $J^*$  is \emph{associated with a block-rook placement} $\Db$
if
$$
\la_\Db = \sum_{(i,j)\in\Db} \la_{ij} E^*_{ij},
$$
moreover, for each block-rook $A_k\times A_\ell\subset \Db$ the submatrix $$\{\la_{ij}:~ i\in A_k,~ j\in A_\ell\}$$ is non-degenerate.\\
\\
\Prop\Num\label{blockJ} \textit{Let $P$ be a parabolic subgroup in  $\GL(n,\Fq)$ of type $(n_1,\ldots, n_s)$. Then any element from  $J$ is equivalent to some element associated with a block-rook placement.}\\
\\
\Proof~ Let  $x\in J$. We shall say a row in  $x$ is \emph{short} if all its non-zero elements are sited in some single block. Replacing  $x$ with an equivalent element one can consider that all rows of $x$ are short and the row systems  of each block-row and of each block-column are linearly independent.

By induction on  $k$, we shall prove that any element  $x\in J$ is equivalent to the element
obeying the condition: all its block-columns with numbers $\ge k$ are associated with a block-rook placement.  We start from the last  $s$th block-column.
 For any  $r\in R_s$, the adjoint action  $x\mapsto rxr^{-1}$ is indeed
  the right-hand multiplication  by  $r^{-1}$; it changes short rows in the last $s$th column and preserve the other short rows.
 Acting by some       appropriate element $r\in R_s$ we change  the last block-column of  $x$ to be associated with some  block-rook placement $\Db^{(s)}$.

 Assume that the statement is proved for  $k+1$; let us prove it for  $k$. Suppose that  $x$ is an element of  $J$, and all its block-columns with numbers $\geq k+1$  are associated with the  block-rook placement $\Db^{(k+1)}$. Let  $L_1,\ldots,L_m$ be the series of block-rooks of the  $k$th block-row sited in the increasing order of block-columns.
 Let  $l_1,\ldots,l_m$ be the orders of the block-rooks  $L_1,\ldots,L_m$.
 Acting by the appropriate element   of the Weyl group $W_{R_k}$
 one can make the rows of the first block-rook $L_1$ to be the first $l_1$ rows in  $I_k$, the rows of  $L_2$ to be the next $l_2$ rows of  $I_k$ and so on. Take $l_{m+1}=n_k - l_1 - \ldots -l_m$.

  Consider the standard parabolic subgroup  $P_k(l)$ in $R_k$ with block sizes  $l=(l_1,\ldots,l_m,l_{m+1})$.

    Define a subgroup $N_{k+1}$  of  $N$  that consists of elements  whose  nonzero non-diagonal  entries are sited in the  block-rows and block-columns with numbers  $> k$.

    For any  $r\in P_k(l)$, there exists  $a\in N_{k+1}$ such that  the element   $ rxr^{-1}a$ obeys the condition:  all its block-columns with numbers $\ge k+1$ are associated with the block-rook placement $\Db^{(k+1)}$, as well as $x$.
    The transformation $x\mapsto   rxr^{-1}a$ acts on the  $k$th  block-column by multiplication $x\mapsto xr^{-1}$. Choosing the appropriate element $r\in P_k(l)$ one can obtain the equivalent element $ rxr^{-1}a$ whose block-columns with numbers  $\geq k$ are associated with a block-rook placement. $\Box$

    Similar statement is true for  $J^*$.\\
    \\
    \Prop\Num\label{blockJstar} \textit{Any element from  $J^*$ is equivalent to some element associated with a block-rook placement.}\\

The following statement proves the Conjecture \ref{conjone} for parabolic subgroups of small block order.
\\
\\
\Theorem\Num\label{classlow} \textit{Let  $P$ be a parabolic subgroup of  $\GL(n,\Fq)$ with blocks of orders  $\leq 2$.\\
 1) For any  $x\in J$ there exists the  rook placement $D$ in $\De_J$ such that  $x$ is equivalent to   $x_D$.\\
2) Two element  $x_D$ and $x_{D'}$ are equivalent if and only if  $D$ and  $D'$ are conjugated with respect to the action of  $W_R$ on  $\De_J$.}\\
\\
\Proof~ \\
 \emph{Item 1.} Applying  Proposition  \ref{blockJ} we consider that  $x=x_\Db$ for some block-rook placement  $\Db$. If all block-rooks from  $\Db$ have order equal to  1, then the statement can be easily proved using the adjoint action of the diagonal subgroup.
Suppose that  $\Db$ has a block-rook of size $2\times 2$.
It can be included  into a  maximal chain   $L_1,\ldots, L_m$ of  $2\times 2$-size block-rooks of $\Db$
subject the requirement  $\col(L_i)=\row(L_{i+1})$. Let  $\row(L_k)=I_{i_k}$.
Suppose that there are two  $1\times 1$ rooks  in the  $i_1$th block-column
and in the  $i_m$th block-row. The other cases can be treated similarly. Using the adjoint action of  the Weyl group of $R_{i_1}$  we can site   the  $1\times 1$-rooks in the  $i_1$th block-column such that the row number and the column number of one of them are less then the ones of other. Analogically, we can site the rooks in
the $i_m$th block-row.
We denote by $X_1,\ldots, X_m$ the chain of submatrices  of $x$  associated with
$L_1,\ldots, L_m$. See the following picture.

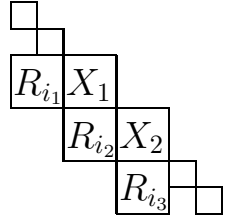
\begin{center}
\begin{picture}(80,80)
\put(0,80){\line(1,0){10}}
\put(0,80){\line(0,-1){10}}
\put(10,80){\line(0,-1){10}}

\put(0,70){\line(1,0){10}}
\put(10,70){\line(1,0){10}}
\put(10,70){\line(0,-1){10}}
\put(20,70){\line(0,-1){10}}
\put(10,60){\line(1,0){10}}

\put(0,60){\line(1,0){20}}
\put(0,60){\line(0,-1){20}}
\put(20,60){\line(0,-1){20}}
\put(0,40){\line(1,0){20}}
\put(1,45){$R_{i_1}$}

\put(20,60){\line(1,0){20}}
\put(20,60){\line(0,-1){20}}
\put(40,60){\line(0,-1){20}}
\put(20,40){\line(1,0){20}}
\put(21,45){$X_1$}

\put(20,40){\line(1,0){20}}
\put(20,40){\line(0,-1){20}}
\put(40,40){\line(0,-1){20}}
\put(20,20){\line(1,0){20}}
\put(21,25){$R_{i_2}$}

\put(40,40){\line(1,0){20}}
\put(40,40){\line(0,-1){20}}
\put(60,40){\line(0,-1){20}}
\put(40,20){\line(1,0){20}}
\put(41,25){$X_2$}

\put(40,20){\line(1,0){20}}
\put(40,20){\line(0,-1){20}}
\put(60,20){\line(0,-1){20}}
\put(40,0){\line(1,0){20}}
\put(41,5){$R_{i_3}$}

\put(60,20){\line(1,0){10}}
\put(60,20){\line(0,-1){10}}
\put(70,20){\line(0,-1){10}}
\put(60,10){\line(1,0){10}}

\put(70,10){\line(1,0){10}}
\put(70,10){\line(0,-1){10}}
\put(80,10){\line(0,-1){10}}
\put(70,0){\line(1,0){10}}

\end{picture}
\end{center}

 As in the proof of Proposition  \ref{blockJ} for any triangular   $2\times 2$-matrix  $t_1\in R_{i_1}$ there exists  $a\in N$ such that   $at_1xt_1^{-1}$ is associated with the same block-rook placement as  $x$.  Similarly, for any triangular   $2\times 2$-matrix  $t_2\in R_{i_m}$ there exists  $b\in N$ such that the element  $t_1xt_1^{-1}b$ is associated with the same block-rook placement as  $x$.
For any chain  $r_{i_1},r_{i_2}, \ldots, r_{i_{m}}, r_{i_{m+1}}$, where $r_{i_{k}}\in R_{i_k}$ and  $r_{i_1}=t_1$,~ $r_{i_{m+1}}=t_2$,
there exists the equivalent element $x'$ with the chain of submatrices   $$X_k'=r_{i_k}X_kr_{i_{k+1}}^{-1}, \quad 1\leq k\leq m.$$
 Then the matrix $Y=X_1\cdots X_m$ is transformed into  $$Y'= X_1'\cdots X_m'= r_{i_1}Yr_{i_{m+1}}^{-1}=t_1Yt_2^{-1}.$$
 There exists triangular matrices   $t_1, t_2$  such that  $Y'$ is a  permutation.
 Choosing consistently $r_{i_2},\ldots, r_{i_{m}}$ we make all $X_1',\ldots, X_{m-1}'$ to be permutations
 (moreover, unit  matrices). Since  $Y'$ is a permutation,  $X_m'$ is also a permutation.

 Proceeding  these transformations for each maximal chain of  $2\times 2$ block-rooks we prove the statement  1).\\
 \emph{Item 2.} Suppose that  $x=x_D$ is equivalent to  $x'=x_{D'}$. Let us show that  $x$ is conjugated to  $x'$ with respect to  $W_R$. For each  $1\leq k<\ell\leq s$ we write $M_{k,\ell}(x)$ for the submatrix of  $x$ with the systems of rows and columns  $I_k\cup\ldots\cup I_\ell$.
 Easy to see that  $x\sim x'$ implies that the submatrices  $M_{k,\ell}(x)$ and $M_{k,\ell}(x')$ are also equivalent. In particular, the ranks of submatrices  $M_{k,\ell}(x)$ and   $M_{k,\ell}(x')$ coincide. Therefore,   $x$ and $x'$ have equal number of ones in each block.

Each nonzero  $2\times 2$-block of the matrix  $x$ can be included into a maximal chain of nonzero  $2\times 2 $-blocks  $X_1,\ldots,X_m$ (as in the previous Item) that can't be extended  due to  one of the following reasons:  all columns with numbers of   $\row(X_1)$ consist of zeros, or all rows with numbers of  $\col(X_m)$ consist of zeroes, or the chain can be extended to
   $$C(x) = \{X_0,X_1,\ldots,X_m,X_{m+1}\},\quad \col(X_i)=\row(X_{i+1}),$$
where the first block  $X_0$ is a nonzero   $\{1\times 2\}$-row and  $X_{m+1}$ is a nonzero   $\{2\times 1\}$-column. We consider the last case; the other cases can be studied similarly.  Let  $M_C(x)$ be the smallest submatrix of $x$ that contains the chain of blocks  $C(x)$; let  $M_C(x')$ be the similar submatrix of  $x'$. Easy to prove that if the matrices  $x$ and $x'$ are equivalent, then  $M_C(x)$ and $M_C(x')$ are conjugated with respect to  $R$.

Suppose that the blocks  $X_1,\ldots,X_m$ are non-degenerate.
Acting on  $x$ and  $x'$ by the appropriate element  $w\in W_R$ we can obtain that $X_0=X_0'=(1,0)$ and the matrices  $X_1,\ldots,X_m, X_1',\ldots, X_m'$ are coincide with the unit matrix. The last blocks  $X_{m+1}$ and $X_{m+1}'$
 can either  $ \left(\begin{array}{c}1\\0\end{array}\right)$ or  $ \left(\begin{array}{c}0\\1\end{array}\right)$.
 If this column is different in $x$ and $x'$, then the matrices  $M_C(x)$ and $M_C(x')$  have different Jordan type, and, therefore, they are not cojugated with respect to  $R$. Hence, the submatrices $M_C(x)$ and $M_C(x')$  are equal (after conjugation by  $w\in W_R$).

 If the chain of blocks $X_1,\ldots,X_m$ contains degenerate matrices, then  we consider the  subchains of non-degenerate blocks in it and deal  with each subchain  as above.

  Considering  each maximal block chain  we conclude that if $x\sim x'$, then  $D$ and  $D'$ are conjugated with respect to $W_R$. $\Box$

The similar statement if true for  $J^*$.\\
\\
\Theorem\Num\label{classlowstar} \textit{ Let  $P$ be a parabolic subgroup in  $\GL(n,\Fq)$ with blocks of orders  $\leq 2$. Then
1) for any  $\la\in J$ there exists a rook placement  $D$ in $\De_J$  such that $\la$ is equivalent to   $\la_D$.\\
2) Two elements  $\la_D$ and $\la_{D'}$ are equivalent if and only if  $D$ and $D'$ are conjugated with respect to the action of  $W_R$ on  $\De_J$.}

\section{Supercharacters}

The statements of this section are proved for an arbitrary parabolic subgroup.  Let $P$ be a parabolic subgroup of $\GL(n,\Fq)$ of type  $(n_1,\ldots, n_s)$.
 Let $\Db$ be a block-rook placement  in $\De_J$, and $\la_\Db$ be an element of  $J^*$ associated with  $\Db$.
The subgroups  $R_\Db$,~ $R_\Db^\circ$ and the characters  $\xi_{\Db,\theta}$,~ $\chi_{\Db,\theta}$ are defined as for   usual rook placement
 $D$.

Notice that  $\RDDr$, $\RDDl$, $R_\Db^\circ$ and $\JDDr$ depend on  $\Db$ and don't depend of the element  $\la_\Db$ associated with the block-rook placement $\Db$; you can see this  in the following example:

\begin{center}
\begin{picture}(460,100)

\put(10,40){\normalsize$\Db=$}
\put(30,100){\line(1,0){100}}
\put(30,80){\line(1,0){100}}
\put(30,60){\line(1,0){100}}
\put(70,40){\line(1,0){60}}
\put(90,0){\line(1,0){40}}
\put(90,20){\line(1,0){40}}

\put(30,100){\line(0,-1){40}}
\put(50,100){\line(0,-1){40}}
\put(70,100){\line(0,-1){60}}
\put(90,100){\line(0,-1){100}}
\put(110,100){\line(0,-1){100}}
\put(130,100){\line(0,-1){100}}

\put(112,63){\small$\otimes$}
\put(120,63){\small$\otimes$}
\put(112,73){\small$\otimes$}
\put(120,73){\small$\otimes$}

\put(132,40){,}

\put(140,40){\normalsize$\JDDr=$}
\put(180,100){\line(1,0){100}}
\put(180,80){\line(1,0){100}}
\put(180,60){\line(1,0){100}}
\put(220,40){\line(1,0){60}}
\put(240,0){\line(1,0){40}}
\put(240,20){\line(1,0){40}}

\put(180,100){\line(0,-1){40}}
\put(200,100){\line(0,-1){40}}
\put(220,100){\line(0,-1){60}}
\put(240,100){\line(0,-1){100}}
\put(260,100){\line(0,-1){100}}
\put(280,100){\line(0,-1){100}}

\multiput(222,82)(10,0){6}%
{\small{*}}
\multiput(222,91)(10,0){6}%
{\small{*}}

\multiput(242,62)(10,0){4}%
{\small{*}}
\multiput(242,71)(10,0){4}%
{\small{*}}

\multiput(242,42)(10,0){4}%
{\small{*}}
\multiput(242,51)(10,0){4}%
{\small{*}}

\put(285,40){,}

\put(300,40){\normalsize$R_\Db^\circ=$}
\put(330,100){\line(1,0){100}}
\put(330,80){\line(1,0){100}}
\put(330,60){\line(1,0){100}}
\put(370,40){\line(1,0){60}}
\put(390,0){\line(1,0){40}}
\put(390,20){\line(1,0){40}}

\put(330,100){\line(0,-1){40}}
\put(350,100){\line(0,-1){40}}
\put(370,100){\line(0,-1){60}}
\put(390,100){\line(0,-1){100}}
\put(410,100){\line(0,-1){100}}
\put(430,100){\line(0,-1){100}}

\multiput(332,82)(10,0){4}%
{\small{$\times$}}
\multiput(332,91)(10,0){4}%
{\small{$\times$}}

\multiput(332,62)(10,0){2}%
{\small{$0$}}
\multiput(332,71)(10,0){2}%
{\small{$0$}}

\put(352,71){\small{$1$}}
\put(362,71){\small{$0$}}
\put(352,62){\small{$0$}}
\put(362,62){\small{$1$}}

\multiput(372,42)(10,0){2}%
{\small{$\times$}}
\multiput(372,51)(10,0){2}%
{\small{$\times$}}

\multiput(392,22)(10,0){2}%
{\small{$\times$}}
\multiput(392,31)(10,0){2}%
{\small{$\times$}}

\multiput(392,2)(10,0){2}%
{\small{$\times$}}
\multiput(392,11)(10,0){2}%
{\small{$\times$}}

\put(412,11){\small{$1$}}
\put(422,11){\small{$0$}}
\put(412,2){\small{$0$}}
\put(422,2){\small{$1$}}

\put(412,31){\small{$0$}}
\put(422,31){\small{$0$}}
\put(412,22){\small{$0$}}
\put(422,22){\small{$0$}}

\end{picture}
\end{center}

Observe that the stabilizer  $R_{\la_\Db}=\{g\in P:~ g\la_\Db g^{-1}=\la_\Db\}$ depends on the choice of  $\la_\Db$ associated with  $\Db$. The same is true for the subgroup  $S\supset R_{\la_\Db} $ defined below.

Let us describe the subalgebra  $\JDDr$.
We say a root  $\al\in \De_J$ is \emph{subordinated } to the root  $\gamma\in \Db$ if they belong to a common row and the block-column number of
           $\al$ is less then the one of  $\gamma$. In this case, if $\al=(i,j)$,~ $\gamma=(i,k)$ with $j<k$, then the positive root  $\beta=(j,k)$ belongs to  $\De_J$.

We say a root  $\al\in \De_J$ is \emph{subordinated } to   $\Db$ if it is subordinated to some root $\gamma\in \Db$.
Denote $$J_{k,\ell}=\mathrm{span}\{E_\al: ~\al\in I_k\times I_\ell\}$$ and  $ \left(J_{\Db, \rt}\right)_{k,l}$ is the intersection of $\JDDr$ with $J_{k,\ell}$.
  Easy to prove the following lemma.
  \\
  \\
  \Lemma\Num ~~ \textit{1) The subalgebra  $J_{\Db, \rt}$ is spanned by the system  $\{E_\al\}$, where $\al\in \De_J$ doesn't subordinate to   $\Db$.\\
  2)    The subalgebra   $\JDDr$ decomposes into a sum of subspaces
   $$ J_{\Db, \rt} = \sum_{k,\ell=1}^m \left(J_{\Db, \rt}\right)_{k,\ell}.$$}\\

Let $S$ be subgroup of $P$ that consists of all  ${g_0}\in P$ such that the action  $\Ad_{g_0}(x)$ preserve $\JDDr$  and stabilize the restriction of $\la_\Db$ on $\JDDr$.
\\
\\
\Lemma\Num\label{Sact}~  \textit{Let ${g_0}\in S$ and ${g_0}=bt$, where $b\in N$  and $t\in R$.
Then \\
1) for any  $r\in \RDDo $ the element $\Ad_t(r)$ also belongs to  $\RDDo$;\\
2) the formula $ \Ad_{g_0}(r) \bmod N = trt^{-1}$ defines the action of  subgroup
 $S$ on  $\RDDo$.}\\
\\
\Proof~ If $r\in \RDDr$, then  $\laDD (r) x = \laDD(x)$ and $$\laDD (\Ad_{g_0}(r) x) = \laDD(x)$$ for any  $x\in \JDDr$.
The element $\Ad_{g_0}(r)$ is presented in the form
\begin{equation}\label{trtr}
\Ad_{g_0}(r)={g_0}r{g_0}^{-1}=btrt^{-1}b^{-1}=r_1b_1,
\end{equation}
where $r_1=trt^{-1}\in R$ and $ b_1=r_1^{-1}br_1b^{-1}$.
Then
\begin{equation}\label{rbrb}
\laDD (r_1b_1 x) = \laDD(x)
\end{equation}
for any  $x\in \JDDr$.
 For the linear form  $\la_\Db$, which is associated with the block-rook placement, the equality (\ref{rbrb}) is equivalent to the pair of equalities
  \begin{equation}\label{rbrb}
\laDD (r_1 x) = \laDD(x),\quad \mbox{and}\quad
\laDD (b_1 x) = \laDD(x),
\end{equation}
 for any  $x\in \JDDr$.

Since $r_1$ belongs to the reductive part  $R$, the first equality is equivalent to $\laDD (r_1 x) = \laDD(x)$ for any  $x\in J$, i.e.  $r\in \RDDr$.
Similarly for  $r\in \RDDl$. Therefore, if  $r\in \RDDo$, then $trt^{-1}\in \RDDo$.  This proves  1) and 2). $\Box$
\\
\\
\Prop\Num\label{sss} ~\textit{Let  $\theta $ and $\theta'$ be  disjoint  $S$-invariant characters of the subgroup  $\RDDo$.
Then the characters  $\chi_{\Db,\theta}$  and $\chi_{\Db,\theta'}$ are disjoint. }\\
\\
\Proof~  Denote  $\xi=\xi_{\Db,\theta}$ and  $\xi' =\xi_{\Db,\theta'}$.
Respectively,  $\chi= \chi_{\Db,\theta} =\Ind(\xi, P_\Db, P)$ and $\chi'= \chi_{\Db,\theta'} =\Ind(\xi', P_\Db, P)$.

It follows from  the Intertwining Number Theorem \cite[теорема 44.5]{CR}
that the characters  $\chi$ and $\chi'$ are disjoint if and only if for any  $g_0\in P$ there exists a subgroup  $H$ in $P_\Db$ such that $ g_0Hg_0^{-1}\subset P_\Db$ and the characters
\begin{equation}\label{state}
\xi^{(g_0)}(h)=\xi(g_0hg_0^{-1})\quad\mbox{and}\quad \xi'(h)\quad \mbox{are~disjoint~on~the~subgroup~} H.
\end{equation}

Order the set of pairs $(k,\ell)$, ~ $1\leq k<\ell\leq m$, as follows
$$(1,m) < (2,m) < \ldots < (m-1,m) < (1,m-1) < (2,m-1) < \ldots < (1,2).$$

The set of roots  $\De_J$ contains the chain of subsets
\begin{equation}\label{dekl}
\De_J^{k,\ell}=\bigcup_{(k_1,\ell_1) < (k,\ell)} I_{k_1}\times I_{\ell_1}.
\end{equation}
The algebra  $J$ contains the chain of ideals   $J^{k,\ell}=\mathrm{span}\{E_\al:~\al\in \De_J^{k,\ell}\}$.
Intersecting  $J^{k,\ell}$ with $\JDDr$ we obtain the chain of subalgebras
\begin{equation}\label{subalg}
J_{\Db,\rt}^{k,\ell}=\JDDr\cap J^{k,\ell}.
\end{equation}
We are aiming to show that one of the following subgroups  $1+J_{\Db,\rt}^{k,\ell}$ or $P_\Db=\RDDo\NDDr$ can be taken as the  subgroup $H$  from (\ref{state}).

  Let $g_0\in P$. Present $g_0$  in the form е $g_0=bt$, where $b\in N$ and $t\in R $.  In order to prove  (\ref{state}), we consider the element  $b$  to be of the form  $b=1+u$, where
   \begin{equation}\label{uuu}
    u=\sum_{\al~\mbox{\small{subordinates}}~\Db} u_\al E_\al,\quad
  u_\al\in\Fq.
   \end{equation}

\emph{Item  1.} Consider the first subalgebra   $J_{\Db,\rt}^{1,m}$ in the chain. The subalgebra  $J_{\Db,\rt}^{1,m}$is an ideal in  $J$; therefore, it is  $\Ad_{g_0}$-invariant.  If the restriction of  $\la$ on $J_{\Db,\rt}^{1,m}$ is not invariant with respect to  $\Ad_{g_0}$, then the disjoint condition  (\ref{state}) is relized with  $H= 1+J_{\Db,\rt}^{1,m}$. This proves that the characters  $\chi$ and $\chi'$ are disjoint.

\emph{Item 2.}  Consider the subalgebra $M=\JDDrkl$ and its previous subalgebra, denote  $M'$, in the chain   (\ref{subalg}).
We simplify the notations:  $\Db_M = \Db\cap\De^{k,\ell}$ and $\Db_{M'}$ is the intersection of $\Db$ with the previous subset in  (\ref{dekl}).
Suppose that the subalgebra $M'$ and the restriction $\la\vert_{M'}$ of the linear form  $\la_\Db$ on  $M'$ are invariant with respect to  $\Ad_{g_0}$.

\emph{Item 2a.} Let us show that the subalgebra  $M$ is invariant with respect to  $\Ad_t$.
 As $t\in R$, the element  $t$ is a block-diagonal matrix $$t=(t_1,\ldots, t_k,\ldots t_\ell,\ldots,t_s).$$
It is sufficient to prove that the  $(k,\ell)$-block  $M_{k,\ell}$ in $M$ is invariant with respect to transformation
$x\to t_kxt_\ell^{-1}$. Since $\JDDr$ is invariant with respect to right-hand multiplication by any elements from  $P$, it is sufficient to show that the block $M_{k,\ell}$ is invariant by with respect to
$x\to t_kx$.

In order to simplify notations, we assume that the $k$th block-row is divided into two  strips of sizes $n_k=n_k'+n_k''$;
 in the first  $n_k'$ rows there is no roots from  $\Db_{M'}$,  and each of the last   $n_k''$ rows  has a root from $\Db_{M'}$.

Then  $M_{k,\ell}$ and the right stabilizer  $R_{\Db_{M'}, \rt}$ consist of the block matrices of the form
 \begin{equation}\label{ttkk}
 M_{k,\ell} = \left(\begin{array}{cc}   * & * \\
  0 & 0 \\
  \end{array}\right),\quad\quad  R_{\Db_{M'}, \rt} =  \left(\begin{array}{cc}   * & * \\
   0 & E \\
 \end{array}\right).
 \end{equation}

 As the subalgebra  $M'$ and the restriction  $\la\vert_{M'}$  are invariant with respect to  $\Ad_{g_0}$,
 the Lemma  \ref{Sact} implies  that  $r\to trt^{-1}$ transform  $R_{\Db_{M'}, \rt}$ into itself.
Then   $t_k$    has the form

   \begin{equation}\label{ttt}
t_k
 =  \left(\begin{array}{cc}   * & * \\
  0 & *\\
 \end{array}\right). \end{equation}
  Therefore,  $M_{k,\ell}$ is invariant with respect to
$x\to t_kx$.

\emph{Item 2b.} Let us show that  $M$ is invariant with respect to  $\Ad_b$, where $g_0=bt$,~  $b=1+u$ and $u$ has the form (\ref{uuu}).

Assume the contrary; then there exists  $E_{jm}\in M_{k,\ell}$ such that $$\Ad_b(E_{jm})\notin M.$$
Since  $M$ is a right ideal in $J$, this is equivalent to  $bE_{jm}\notin M$. By (\ref{uuu}), we obtain $$bE_{jm}= \sum_{(i,j)~\mbox{\small{subordinate}}~\Db}   u_{ij}E_{ij}E_{jm} = \sum_{(i,j)~\mbox{\small{subordinate}}~\Db}u_{ij}E_{im}\notin M.$$

There exists a root  $\al=(i,j)$ such that   $u_{ij}\ne 0$ and it subordinates to some root $\gamma=(i,p)\in\Db_{M'}$, where $m<p$, together with  $(i,m)$.
As   $E_{jm}\in M_{k,\ell}$, the root element $E_{jp}$  belong to  $M$; more precisely it belongs to some block $ M_{k,\ell'}$, where
$\ell<\ell'$.
 Take $Y=\Ad_{t}^{-1}(E_{jp})$. As well as $E_{jp}$, the element $Y$ belongs to the first  $n_k'$ rows of the  $k$th block-row.
Therefore,  $Y\in M_{k,\ell'}$ and   $\la_\Db(Y)=0$.
According to the assumption of Item  2  the subalgebra  $M'$ and the linear form  $\la\vert_{M'}$ are invariant with respect to   $\Ad_{g_0}$; then
$$Ad_{g_0}(Y)=\Ad_b(E_{jp})\in M' \quad \mbox{and}\quad \la_\Db(\Ad_{g_0}(Y))=\la_\Db(Y)=0.$$

 The root $\gamma$ belongs to some block-rook  $A\times B\in \Db_{M'}$. Therefore, the root $\al$ can be included into a system of roots  $\al_1=(i_1,j), \ldots, \al_{|A|}=(i_{|A|},j)$ being subordinated to to the same block-rook  $A\times B$, where $A=\{i_1,\ldots,i_{|A|}\}$.
Similarly to  $\al$ we have $E_{i_cj}E_{jp}=E_{i_cp}$ and $(i_c,p)\in A\times B$ for any  $1\leq c\leq |A|$.

The element $Ad_{g_0}(Y)$ can be presented in the form
$$\Ad_{g_0}(Y)= \Ad_b(E_{jp})= \sum_{1\leq c\leq |A|} u_{\al_c} E_{i_cp} + Z,$$
where $Z\in M'$ and $\la_\Db(Z)=0$.

Then $$\la_\Db(\Ad_{g_0}(Y))= \sum _{1\leq c\leq |A|} u_{\al_c} \la_{i_cp} =0,$$
where $\la_{i_cp}=\la_\Db(E_{i_cp})$.
Since the matrix  $\{\la_{i_cp}: ~ i_c \in A, ~p\in B\}$ is non-degenerate (see Definition  \ref{assblock}), we have
$ u_{\al_c}=0$ for any  $1\leq c\leq |A|$. This contradicts to choice of the root  $\al$  with  $u_\al\ne 0$.
We conclude that  $M$ is invariant with respect to  $\Ad_b$.
The statement of Item 2b) is proved.

\emph{Item 2c}. So, the subalgebra  $M$ is invariant with respect to  $\Ad_t$ and $\Ad_b$. Therefore, it is  $\Ad_{g_0}$-invariant. If  $\Ad_{g_0}(\la_\Db\vert_M)\ne \la_\Db\vert_M$, then the characters $\Ad_{g_0}\xi$ and $\xi'$  don't coincide on the subgroup  $1+M$; the condition  (\ref{state}) is fulfilled. The characters  $\chi$ are $\chi'$ disjoint.

\emph{Item 3}.
 Suppose that  $\Ad_{g_0}$ preserve all subalgebras $J_{\Db,\rt}^{k,\ell}$ and the linear form  $\la_\Db$ on them.
 Let us show that in this case the subgroup $H = P_\Db$ is invariant with respect to $\Ad_{g_0}$ and
  the disjoint condition  (\ref{state}) is valid on $H$.

\emph{Item 3a}. Suppose as above that  $b=1+u$, and $u$ has the form  (\ref{uuu}). Let us prove that, for any  $r\in \RDDr$, the element $r^{-1}ur$ is presented in the form  $r^{-1}ur=u+v$, where  $v\in \JDDr$ and $\la_\Db(v)=0$.

It is sufficient to show that  the equality $r_k^{-1}E_\al r_\ell = E_\al+v$, where  $v\in \JDDr$ and $\la_\Db(v)=0$  is valid for any  $(k,\ell)$ and for any  $E_\al\in J_{k,\ell}$ with $u_\al\ne 0$ (see (\ref{uuu})).

Similarly to  Item 2a, in order to simplify notations, we suppose that the $k$th row is decomposed into tow strips  $n_k=n_k'+n_k''$ such that there is no roots from $\Db$ in the first  $n_k'$ rows of the $k$th block-row,  and  in each of  the others  $n_k''$ row there is a root from $\Db$. Analogically, we decompose the   $\ell$th block-row $ n_\ell=n_\ell'+n_\ell''$.
Since  $r\in \RDDr$, we have
$$r_{k}= \left(\begin{array}{cc} A_1&B_1\\0&E_{n_k''}\end{array}\right), \quad r_{\ell} = \left(\begin{array}{cc} A_2&B_2\\0&E_{n_\ell''}\end{array}\right),$$
where $E_{n_k''}$ and $E_{n_\ell''}$ are the unit matrices of orders  $n_k''$ and $n_\ell''$.

As the root $\al=(i,j)$ is subordinated to  $\Db$, the row $i$ is one of the last  $n_k''$ rows of the $k$th block-row. If the  $j$th row is one of the first  $n_\ell'$ rows of the  $\ell$th block row, then there exists a root   $(j,p)$ such that  $E_{jp}\in \JDDr$ and $(i,p)\in\Db$. Arguing as in the Item 2b, we may show that the assumption of Item  3 contradicts to  $u_\al\ne 0$.
Therefore, the   $j$th row is one of the last  $n_\ell''$ rows of the  $\ell$th block row, and
 $E_\al$ (as an element of the block  $J_{k,\ell}$) has the form $E_\al = \left(\begin{array}{cc} 0&0\\0&U\end{array}\right)$,
 where  $U$ is a submatrix of size $n_k''\times n_\ell''$.
 We have
$$
 r_{k}^{-1} E_\al  r_{\ell} = \left(\begin{array}{cc} A_1&B_1\\0&E_k\end{array}\right)^{-1} \left(\begin{array}{cc} 0&0\\0&U\end{array}\right) \left(\begin{array}{cc} A_2&B_2\\0&E_\ell\end{array}\right) = \left(\begin{array}{cc} 0 & *\\0&U\end{array}\right) = E_\al + \left(\begin{array}{cc} 0& *\\0&0\end{array}\right).$$
This proves the statemnet of Item  3a.

\emph{Item 3b}.  From  (\ref{trtr}) we obtain
$g_0rg_0^{-1}= r_1b_1$, where $r_1=trt^{-1}\in R$ and $b_1=r_1^{-1}br_1b^{-1}\in N$. By the condition of Item  3, the element $g_0$ belongs to  $S$. According to Lemma \ref{Sact} the element $r_1$  belongs to $\RDDo$. Let us show that $b_1\in \NDDr$ and $\xi(b_1)=1$.

Really, from Item 3a, we have  $r_1^{-1}br_1=1+r_1^{-1}ur_1 =1+ u+v$, where  $v\in \JDDr$ and $\la_\Db(v)=0$.
Then
$$b_1=(1+u+v)(1+u)^{-1})=1+v(1+u)^{-1}=1+v+vu_1$$
for some  $u_1\in J$.
Since $\JDDr$ is a right ideal of $J$, we have $vu_1\in  \JDDr$ and $\la_\Db(vu_1)=0$. Hence $b_1\in \NDDr$. Then
$$\xi(b_1)=\eps^{\la_\Db(v+vu_1)}= \eps^{\la_\Db(v)}\eps^{\la_\Db(vu_1)} = 1.$$

\emph{Item 3c}. From Item 3b) $g_0rg_0^{-1}\in P_\Db$. Let us prove the disjoint condition (\ref{state}).
For $h=ra\in P_\Db$, where $r\in R_\Db^0$ and $a\in \NDDr$, we obtain
 $$\xi^{(g_0)}(h) = \xi(g_0rag_0^{-1}) = \xi(g_0rg_0^{-1}g_0ag_0^{-1})=\theta(trt^{-1})\xi(b_1)\xi(g_0ag_0^{-1}).$$
By condition, the character  $\theta$  is invariant with respect to  $S$; then $\theta(trt^{-1})=\theta(r)$. By Item 3b, $\xi(b_1)=1$.
By the condition of  Item 3,   $\xi(g_0ag_0^{-1})=\xi(a)$. Finally, we get $$\xi^{(g_0)}(h)=\theta(r)\eps^{\la_\Db(a-1)}$$.

 Comparing with  $\xi'(h)=\theta'(r)\eps^{\la_\Db(a-1)}$, we conclude that the disjoint condition  (\ref{state}) is valid for  $H=P_\Db$.
The characters  $\chi$ and $\chi'$ are disjoint. $\Box$
\\
\\
\Prop\Num\label{supersuper} \textit{The characters   $\{\chi_{D,\theta}\} $ are constant on the equivalence classes $\{K(g)\}$.}\\
\\
\Proof~ Let  $\chi=\chi_{\Db,\theta}$ and $\xi=\xi_{\Db,\theta}$.

\emph{Item 1}. Let $g\in P_\Db$, ~$a\in N$. Let us prove that if   $g'=1+(g-1)a\in P_\Db$, then $\chi(g')=\chi(g)$.

As $g\in P_\Db$, we get $g=r+x$, where $r\in \RDDo$ and $x\in \JDDr$. Let $a=1+u$, where $u\in J$.
 Then $g'= r+y$, where $r\in \RDDo$ and $y=(r-1)u + x+ xu\in J$.  The equalities  $$\la_\Db r =\la_\Db,\quad \la_\Db(r-1)u=0, \quad \la_\Db x=0~~\mbox{and}~~ \la_\Db xu=0$$ imply
 $\la_\Db y = \la_\Db (r-1)u +\la_\Db x + \la_\Db xu = 0$. Therefore,  $y\in \JDDr$ and $g'\in P_\Db$.

We obtain
\begin{equation}\label{xixi}
    \xi(g') = \theta(h)\varepsilon^{\la_\Db((r-1)u+x+xu)}=\theta(r)\varepsilon^{\la_\Db(x)}=\xi(g).
\end{equation}
Denote $\Lambda(g) = \{s\in G|~ s^{-1}gs\in P_\Db\}$.

Let us show that  $\Lambda(g') = \Lambda(g)$ for any  $g\in P_\Db$. Indeed,  $s^{-1}g's = 1+ (s^{-1}gs-1)s^{-1}as$.
As we saw above  if  $ s^{-1}gs\in P_\Db$, then $s^{-1}g's\in P_\Db$. This proves  $\Lambda(g)\subset\Lambda(g')$. The converse  inclusion can be proved analogically.

The formula  (\ref{xixi}) implies
\begin{equation}
\xi(s^{-1}g's) = \xi(s^{-1}g s),
\end{equation}
for any $s\in \Lambda(g)$.
Hence
\begin{equation}\label{chichi}
\chi(g') = \frac{1}{|P_\Db|}\sum_{s\in \Lambda(g')} \xi(s^{-1}g's)= \frac{1}{|P_\Db|}\sum_{s\in \Lambda(g)} \xi(s^{-1}gs) = \chi(g).
\end{equation}
\emph{Item 2}. The characters are constant on conjugacy  classes. Therefore, for any $g_0\in P$, we have
$\chi(g_0gg_0^{-1})=\chi(g)$.
Then $\chi $ is constant on  $K(g)$ for any  $g\in P_\Db$.
If an equivalence class has empty intersection with  $P_\Db$, then the value of $\chi$   is zero on this class. $\Box$

\section{The case of blocks of small orders}

Let $P$ be a parabolic subgroup with blocks  $(n_1,\ldots, n_s)$, where $n_k \leq 2$.
We prove conjectures of section 3 in this case.

Let  $Cl_2$ stand for the set of representatives of conjugasy classes for the group $\mathrm{GL}(2,\Fq)$
that contains the matrices of the form а $\left(\begin{array}{cc} a&0\\0&1\end{array}\right)$ with $a\ne 0$ and  $\left(\begin{array}{cc} 1&1\\0&1\end{array}\right)$.\\
\\
\Notation\Num\label{defb}
 Denote by  $\Bc$ the set of pairs  $(D,\rho)$, where $D$ is a representative of the $W_R$-orbits
  on rook placements in  $\De_J$, and  $\rho=(\rho_1,\ldots,\rho_s)$ is an element of  $R$ such that  each block  $\rho_*$ obeys one of the following conditions:\\
1) The order of the block  $\rho_*=(a_{ii})$ equals to  1.  \\
 1a)~  If there are no roots of $D$ in the   $i$th row and in the  $i$th column, then  $a_{ii}$ is an arbitrary element of  $\Fq^*$.\\
  1b)~ If there is a root of $D$ in the   $i$th row or in the  $i$th column, then  $a_{ii}=1$.\\
  2) The order of the block $\rho_*=\left(\begin{array}{cc} a_{i,i}&a_{i,i+1}\\a_{i+1,i}&a_{i+1,i+1}\end{array}\right)$ equals to 2. \\
 2a)~ If there are no roots of $D$ in the  $\rho_*$ block-row and in the $\rho_*$ block column,
 then $\rho_*$ is an arbitrary element of  $Cl_2$.\\
  2b)~ If there is a root of $D$ in the $\{i+1\}$th row, and
   there is  no roots  of $D$ in the $i$th row and in the  $\rho_*$ block-column, then   $\rho_* = \left(\begin{array}{cc} a&0\\0&1\end{array}\right)$ or  $\rho_* = \left(\begin{array}{cc} 1&1\\0&1\end{array}\right)$.\\
 If we invert here  $i$ and $i+1$,  then  we  substitute  $\rho_*$ for  $\sigma \rho_* \sigma^{-1}$, where
$\sigma=\left(\begin{array}{cc} 0&1\\1&0\end{array}\right)$.\\
  2c)~ If there is a root of $D$ in the $i$th column, and
   there is  no roots  of $D$ in the $\{i+1\}$th column and in the  $\rho_*$ block-row, then   $\rho_* = \left(\begin{array}{cc} 1&0\\0&a\end{array}\right)$ or  $\rho_* = \left(\begin{array}{cc} 1&1\\0&1\end{array}\right)$.\\
 If we invert here  $i$ and $i+1$,  then  we  substitute  $\rho_*$ for  $\sigma \rho_* \sigma^{-1}$.\\
  2d) If the  $\{i+1\}$th row and the
  $\{i+1\}$th column have roots of  $D$, then   $\rho_* = \left(\begin{array}{cc} a&0\\0&1\end{array}\right)$.
  If we substitute here $i+1$ for $i$, then $\rho_*$ is substituted for  $\sigma \rho_* \sigma^{-1}$.\\
  2e)~ If the   $\{i+1\}$th row and the
  $i$th column have roots of  $D$, then $\rho_* = \left(\begin{array}{cc} 1&1\\0&1\end{array}\right)$ or $\rho_* = \left(\begin{array}{cc} 1&0\\0&1\end{array}\right)$.\\
  If  we invert here $i+1$ and $i$, then  $\rho_*$ is  substituted for  $\sigma \rho_* \sigma^{-1}$.\\
  2f)~ If there are two roots of $D$ either in the $\rho_*$ block-row, or in the $\rho_*$ block-column, then $\rho_* = \left(\begin{array}{cc} 1&0\\0&1\end{array}\right)$.\\

By any pair $(D,\rho)\in \Bc$, we  define an element  $g_{D,\rho}=\rho+x_D$, where $x_D$ from (\ref{xD}).
 Here are some example of matrices of the form  $g_{D,\rho}$:

\begin{center}
\begin{picture}(460,40)
\put(0,40){\line(1,0){40}}
\put(0,30){\line(1,0){40}}
\put(10,10){\line(1,0){30}}
\put(30,0){\line(1,0){10}}

\put(40,0){\line(0,1){40}}
\put(30,0){\line(0,1){40}}
\put(10,10){\line(0,1){30}}
\put(0,30){\line(0,1){10}}

\put(2,32){\small{$1$}}
\put(32,2){\small{$1$}}
\put(32,32){\small{$1$}}
\put(15,15){{$A$}}

\put(60,15){{,}}


\put(80,40){\line(1,0){40}}
\put(80,30){\line(1,0){40}}
\put(90,10){\line(1,0){30}}
\put(110,0){\line(1,0){10}}

\put(120,0){\line(0,1){40}}
\put(110,0){\line(0,1){40}}
\put(90,10){\line(0,1){30}}
\put(80,30){\line(0,1){10}}

\put(82,32){\small{$1$}}
\put(112,2){\small{$1$}}

\put(102,32){\small{$1$}}
\put(112,12){\small{$1$}}

\put(92,22){\small{$a$}}
\put(102,12){\small{$1$}}

\put(90,20){\line(1,0){30}}
\put(100,10){\line(0,1){30}}

\put(140,15){{,}}


\put(160,40){\line(1,0){40}}
\put(160,30){\line(1,0){40}}
\put(170,10){\line(1,0){30}}
\put(190,0){\line(1,0){10}}

\put(200,0){\line(0,1){40}}
\put(190,0){\line(0,1){40}}
\put(170,10){\line(0,1){30}}
\put(160,30){\line(0,1){10}}

\put(162,32){\small{$1$}}
\put(192,2){\small{$1$}}

\put(172,32){\small{$1$}}

\put(172,22){\small{$1$}}
\put(182,12){\small{$a$}}

\put(170,20){\line(1,0){30}}
\put(180,10){\line(0,1){30}}

\put(220,15){{,}}


\put(240,40){\line(1,0){40}}
\put(240,30){\line(1,0){40}}
\put(250,10){\line(1,0){30}}
\put(270,0){\line(1,0){10}}

\put(280,0){\line(0,1){40}}
\put(270,0){\line(0,1){40}}
\put(250,10){\line(0,1){30}}
\put(240,30){\line(0,1){10}}

\put(242,32){\small{$1$}}
\put(272,2){\small{$1$}}


\put(252,22){\small{$1$}}
\put(262,12){\small{$1$}}

\put(250,20){\line(1,0){30}}
\put(260,10){\line(0,1){30}}

\put(252,32){\small{$1$}}
\put(262,22){\small{$1$}}
\put(272,12){\small{$1$}}

\put(300,15){{,}}


\put(320,40){\line(1,0){40}}
\put(320,30){\line(1,0){40}}
\put(330,10){\line(1,0){30}}
\put(350,0){\line(1,0){10}}

\put(360,0){\line(0,1){40}}
\put(350,0){\line(0,1){40}}
\put(330,10){\line(0,1){30}}
\put(320,30){\line(0,1){10}}

\put(322,32){\small{$1$}}
\put(352,2){\small{$1$}}


\put(332,22){\small{$1$}}
\put(342,12){\small{$1$}}

\put(330,20){\line(1,0){30}}
\put(340,10){\line(0,1){30}}

\put(332,12){\small{$1$}}
\put(352,22){\small{$1$}}

\put(380,15){{,}}


\put(400,40){\line(1,0){40}}
\put(400,30){\line(1,0){40}}
\put(410,10){\line(1,0){30}}
\put(430,0){\line(1,0){10}}

\put(440,0){\line(0,1){40}}
\put(430,0){\line(0,1){40}}
\put(410,10){\line(0,1){30}}
\put(400,30){\line(0,1){10}}

\put(402,32){\small{$a$}}
\put(432,2){\small{$b$}}
\put(415,15){{$A$}}

\put(450,15){{.}}
\end{picture}
\end{center}
\Theorem\Num\label{theosupercl}  \textit{Let $P$ be a parabolic subgroup with blocks of orders  $\leq 2$.
Then\\
1) an arbitrary element  $g$ of $P$ is equivalent (see  definition \ref{superclass}) to some element  $g_{D,\rho}$ constructed by  the pair  $(D,\rho)\in \Bc$;\\
2) the element  $g_{D,\rho}$ is uniquely defined by  $g$;\\
3) the group $P$ decomposes into the system of subsets $K_{D,\rho}= K(g_{D,\rho})$.}\\
\\
\Proof~ 1) Let $g=r+x$, where $r\in R$,~ $x\in J$. Acting by  $\Ad_R$ we obtain  $r$ of the form $r=(r_1,\ldots,r_s)$, where each $r_*$ is an element of  $Cl_2$.
There exist  $A,B\in N$  such that the element  $g'=1+A(g-1)B$ satisfies the conditions:\\
i) if the block $r_*$ is such that   $\det(r_*-1)\ne 0$, then  $g'=r+x'$, where $x'$ has  the zero  $r_*$ block-row and zero $r_*$ block-column;\\
ii) if $r_*=\left(\begin{array}{cc} a_{i,i}&a_{i,i+1}\\a_{i+1,i}&a_{i+1,i+1}\end{array}\right)$, where $a_{i,i}=a\notin\{0,1\}$, ~$a_{i,i+1}=a_{i+1,i}=0$ and $a_{i+1,i+1}=1$, then $x'$ has the zero $i$th row and  zero $i$th column;
\\
iii) if $a_{i,i}=a_{i,i+1}=a_{i+1,i+1}=1$, ~$a_{i+1,i}=0$, then $x'$ has the zero $i$th row and zero $\{i+1\}$th column.

Arguing as in the proof of Theorem \ref{classlow} one can show that the element  $g'=r+x'$ is equivalent to   $g''=r+x_D$, where all roots of $D$ belong to rows and columns listed above. This implies  1). The statement  2) is proved similarly to the proof of Item  2 from  Theorem  \ref{classlow}. $\Box$\\
\\
\Notation\Num ~ Denote by  $\Ac$ the set of pairs  $(D,\theta)$, where $D$ is a representative of
   the $W_R$-orbits on the rook placements of $\De_J$, and  $\theta$ is an  $R_D$-irreducible representation of the subgroup  $P_D^\circ$.\\
   \\
\Theorem\Num\label{maintheo}~  \textit{Let $P$ be a parabolic subgroup with blocks of orders  $\leq 2$.
Then the system of characters  $$\{\chi_{D,\theta}:~ (D,\theta)\in \Ac\}$$ and the partition of the subgroup  $P$ into the subsets $$\{K_{D,\rho}:~ (D,\rho)\in \Bc\}$$ give rise to a supercharacter theory of the group  $P$.}\\
\\
\Proof~ \\
\emph{Item 1. } Let us show that  $|\Ac|=|\Bc|$. It is sufficient to prove that  the fibers $\Ac_D$ and $\Bc_D$ have equal number of elements for any $D$. The subgroup  $R_D^\circ$ decomposes into a product $R_D^\circ=R_{D,1}^\circ\times\ldots\times R_{D,s}^\circ$.
We denote by  $\al_*$ and $\beta_*$  a contribution of the block $R_*$ in $|\Ac_D|$ and $|\Bc_D|$.

Let us show that $\al_*=\beta_*$ in each of the cases listed in Definition  \ref{defb}.
If the order of block  $R_*$ equals to  1, then in the case  1a  $\rho_*$ runs through  $\Fq^*$ and $\beta_*=q-1$.
On the other hand, in this case,  $ R^\circ_{D,*}=\Fq^*$ and $\al_*=q-1$. Then $\al_*=\beta_*=q-1$.
In the case  1b, we have $\al_*=\beta_* = 1$.

Let the order of the block $R_*$ equals to  2. Then in the case 2a, we obtain  $\al_* = \beta_*$ since  $R_{D,*}=\mathrm{GL}(2,\Fq)$ and the number of conjugacy classes equals to the number of irreducible characters. In the case 2b, the subgroup  $R^\circ_{D,*}$ coincides with the subgroup of matrices  $\left\{\left(\begin{array}{cc} a&b\\0&1\end{array}\right)\right\}$.
The projection   $R_{D,*}$ of the stabilizer  $R_D$ onto $R_*$ is contained in the normalizer
\begin{equation}\label{norm}
\mathrm{Norm}_{D,*} = \left\{\left(\begin{array}{cc} a&b\\0&c\end{array}\right)\right\}
\end{equation}
of the subgroup $R^\circ_{D,*}$ in $\GL(2,\Fq)$.
Each irreducible character of   $R^\circ_{D,*}$ is invariant with respect to  $\mathrm{Norm}_{D,*}$ and, therefore, with respect to  $R_{D,*}$. Then  $\al_*=\beta_* = q$.

The cases  2b, 2c, 2d are processed analogically. Consider the case 2e. Here  $\beta_{*}=2$.
The subgroup  $R^\circ_{D,*}$ coincides with the subgroup of matrices $\left\{\left(\begin{array}{cc} 1&b\\0&1\end{array}\right)\right\}$.
The subgroup  $R_D$ contains the subgroup  $H_i$ of the diagonal matrices $a_{11}=a,\ldots, a_{ii}=a,~ a_{i+1,i+1}=b,\ldots, a_{nn}=b.$
 The projection of $H_i$ onto   $R_*$  is the subgroup  $H_{i,*}=\left\{\left(\begin{array}{cc} a&0\\0&b\end{array}\right)\right\}$. Therefore,  $R_{D,*}$ coincides with  $\mathrm{Norm}_{D,*}$ from (\ref{norm}).
The number of $\mathrm{Norm}_{D,*}$-irreducible characters of $R^\circ_{D,*}$ equals to  2. We have
$\al_*=\beta_* = 2$.

Finally, in the last case 2f, we get  $R^\circ_{D,*}=\{1\}$ and $\al_*=\beta_* = 1$. This proves the statement  of Item 1.

\emph{Item 2.} Let us show that  the characters $\{\chi_{D,\theta}\}$ are pairwise disjoint.

Consider the case     $D$ and $D'$ belong to different  $W_R$-orbits. The restriction of the character  $\{\chi_{D,\theta}\}$ on the subgroup $N$ is a sum of characters of the form  $\chi_{wD}$, where $w\in W_R$, and  $\chi_{D}$ is a supercharacter of the subgroup $N$ constructed by $\la_D$ in  the paper \cite{DI} of  P.Diaconis and  I.M.Isaacs. The characters $\chi_D$ and $\chi_{D'}$ are disjoint if they belong to different left-right $N$-orbits.
As $D$ and $D'$ belong to different  $W_R$-orbits, then  $\la_D$ and $\la_{D'}$ are not equivalent (see Theorem \ref{classlowstar}), and
the restrictions of  characters $\{\chi_{D,\theta}\}$ and $\{\chi_{D',\theta'}\}$ on $N$ are disjoint. Then the characters $\{\chi_{D,\theta}\}$ and $\{\chi_{D',\theta'}\}$  are disjoint.

 Suppose that $D=D'$. Apply the Proposition  \ref{sss}. It is sufficient to prove that  the $R_D$-invariant characters  $\theta$ and  $\theta'$ are invariant with respect to the subgroup  $S$
 defined in the previous section. Indeed, let   $S_*$ be the projection of the subgroup  $S$ onto the block  $R_*$. Then  $R_{D,*}\subseteq S_*\subseteq  \mathrm{Norm}_{D,*}$. In each of case  $1a-1b$ and $2a-2f$ the  $R_D$-invariancy of characters  implies their  $S$-invariancy. The characters  $\{\chi_{D,\theta}\}$ and $\{\chi_{D',\theta}\}$ are disjoint.

\emph{Item 3.} By Proposition \ref{supersuper}, the characters  $\{\chi_{D,\theta}\} $ are constant on the equivalence classes
$\{K_{D,\rho}\}$; this proves condition   1) of Definition  \ref{defsuptheory}. Finally,   $\{1\}$ is an equivalence class of the unit element of the subgroup  $P$.

According to the Definition  \ref{defsuptheory} the system of characters  $\{\chi_{D,\theta}\} $ and the partition of $P$ into equivalence classes  $\{K_{D,\rho}\}$ give rise to a supercharacter theory of the subgroup $P$. $\Box$

\end{document}